\documentclass{article}
\usepackage{amsmath,amsfonts,amsthm,amssymb}
\usepackage{graphicx,color}
\numberwithin{equation}{section}
\newtheorem{thm}{Theorem}
\newtheorem{cor}[thm]{Corollary}
\theoremstyle{remark}
\newtheorem{expl}{Example}[section]
\def\R{\mathbb R}
\def\D{\mathrm D}\def\d{\mathrm d}
\def\N{\mathbb N}
\let\<\langle
\let\>\rangle
\DeclareMathOperator*{\slim}{s-lim}
\author{Michael Reissig \and 
Jens Wirth\thanks{This work was partially supported
by the government of the state of Saxony with a Landesgraduiertenstipendium, Grant G 151}
}
\title{$L^p$--$L^q$ decay estimates for 
wave equations with monotone time-dependent dissipation}
\begin{document}
\maketitle
{\centering\it 
Institute of Applied Analysis, 
TU Bergakademie Freiberg, \\
09596 Freiberg, Germany\\}
\begin{abstract}
This expository article is intended to give an overview about recently
achieved results on asymptotic properties of solutions to the Cauchy
problem
\begin{equation*}
  u_{tt}-\Delta u+b(t)u_t =0,\qquad u(0,\cdot)=u_1,\quad \D_tu(0,\cdot)=u_2 
\end{equation*}
for a wave equation with time-dependent dissipation term. The results are 
based on structural properties of the Fourier multipliers representing its 
solution.

The article explains the general philosophy behind the approach.
\end{abstract}

\section{Introduction}
\paragraph{Strichartz and Matsumura type estimates.}
To prove global existence results for small data solutions to Cauchy problems
for nonlinear wave equations so-called Strichartz' decay estimates for the 
energy $||(\nabla u(t,\cdot), u_t(t,\cdot))||_q$ based on the $L^q$-norm, 
$q\geq 2$, are an essential ingredient. For the free wave equation
$$  u_{tt}-\Delta u=0,\qquad u(0,\cdot)=u_1,\quad \D_t u(0,\cdot)=u_2, $$
one obtains, \cite{Str69}, \cite{vW71},
\begin{equation}\label{eq:Strichartz}
||(\nabla u(t,\cdot), u_t(t,\cdot))||_q
\leq C(1+t)^{-\frac{n-1}2\left(\frac1p-\frac1q\right)}
||(\<\D\>u_1,u_2)||_{L^{p,r_p}}
\end{equation}
on the conjugate line $pq=p+q$, $q\in[2,\infty]$ and with 
$r_p>n\left(1/p-1/q\right)$. Here and thereafter we denote by 
$L^{p,r}(\R^n)=\<\D\>^{-r}L^p(\R^n)$ the Bessel potential space
of order $r$ over $L^p(\R^n)$. For $p=q=2$ the estimate is related to
the conservation of the energy while for $p=1$ and $q=\infty$ the uniform
decay of the energy of the solutions follows for space dimension $n>1$.

If we include a further constant dissipation term
$$  u_{tt}-\Delta u+u_t=0,\qquad u(0,\cdot)=u_1,\quad \D_tu(0,\cdot)=u_2, $$
the results of Matsumura, \cite{Mat76}, imply the corresponding estimate
\begin{equation}\label{eq:Matsumura}
||(\nabla u(t,\cdot), u_t(t,\cdot))||_q
\leq C(1+t)^{-\frac{n}2\left(\frac1p-\frac1q\right)-\frac12}
||(\<\D\>u_1,u_2)||_{L^{p,r_p}}
\end{equation}
under the same assumptions on $p$, $q$, $r_p$. The dissipation term 
we feel in the  occurrence of the further decay factor $(1+t)^{-\frac12}$ 
in the $L^2$--$L^2$ estimate and in the different constant $n/2$
in front of $(1/p-1/q)$ instead of $(n-1)/2$. 

It is a natural question to ask for generalisations of these two estimates
to classes of variable-coefficient dissipation terms. 
Under suitable assumptions on the coefficient function so-called weighted 
energy inequalities provide a tool to obtain $L^2$--$L^2$ estimates. They are
not used in our approach, so we refer to \cite{Mat77}, \cite{Ues79} and 
\cite{HN03} and the references cited therein. Weighted energy inequalities
have the disadvantage that they do not give information on structural 
properties of the solution operator to the Cauchy problem.

\paragraph{An intermediate case.} For the special model problem 
\begin{equation}
  u_{tt}-\Delta u+\frac{\mu}{1+t}u_t=0,
\qquad u(0,\cdot)=u_1,\quad \D_t u(0,\cdot)=u_2,
\end{equation}
where $\mu$ is a non-negative constant, we can use a relation to Bessel's
differential equation and obtain an explicit representation of its solution
in terms of special functions, \cite{Wir02}. A careful phase space analysis
yields for it the $L^p$--$L^q$ decay estimate
\begin{multline}\label{eq:Wirth}
||(\nabla u(t,\cdot), u_t(t,\cdot))||_q\\
\leq C(1+t)^{\max\{-\frac{n-1}2\left(\frac1p-\frac1q\right)-\frac\mu2,\,
-n\left(\frac1p-\frac1q\right)-1\}}
||(\<\D\>u_1,u_2)||_{L^{p,r_p}}
\end{multline}
on the conjugate line $pq=p+q$, $q\in[2,\infty]$ and with 
$r_p>n\left(1/p-1/q\right)$. The estimate allows
interesting observations and there arise several questions related to it.
\begin{itemize}
\item If $p=q=2$, then small values of $\mu$ have a direct influence on the 
decay rate, while we do not feel large values of $\mu$. The value $\mu=2$ is
critical for such $L^2$--$L^2$ estimates. 
\item If $\mu\leq2$, estimate \eqref{eq:Wirth} generalises 
\eqref{eq:Strichartz}. The dissipation does not destroy the structure of the 
Strichartz decay estimate, it implies a further decay factor 
$(1+t)^{-\frac\mu2}$. In this case we will call the dissipation 
{\em non-effective}.
\item If $\mu\geq n+3$, then the decay order is $-n(1/p-1/q)-1$. In this case
the structure of the estimate changed completely and we will call
the dissipation {\em effective}, therefore. 
Does there exist a relation between this estimate and \eqref{eq:Matsumura}?
\item For $2<\mu<n+3$ there appears a mixture of both situations. The
structure of the decay rate dependents on the choice of $p$ and $q$.
\end{itemize}

\paragraph{Main objectives and basic assumptions.}
We will extent these observations to a broader class of monotone dissipation
terms $b(t)u_t$. For the coefficient function $b=b(t)\in C^\infty(\R^n)$ 
we impose the following conditions, 
\begin{itemize}
\item[(H1)] $b(t)\geq 0$,
\item[(H2)] $b(t)$ is monotone in $t$,
\item[(H3)] for all $k\in\mathbb N_+$ we have
$$ \left|\frac{\d^k}{\d t^k}b(t)\right|\leq C_k b(t)\left(\frac1{1+t}\right)^k $$
with suitable constants $C_k$.
\end{itemize}

Estimate \eqref{eq:Wirth} turns out to be the intermediate
case in between two different scenarios occurring for the Cauchy 
problem
\begin{equation}\label{eq:CP}
  u_{tt}-\Delta u+b(t)u_t =0,\qquad u(0,\cdot)=u_1,\quad \D_tu(0,\cdot)=u_2 
\end{equation}
with such a time-dependent dissipation term. In the following we will sketch 
the main ideas of the approach together with their consequences on the
$L^p$--$L^q$ decay of solutions and their energy.

The results are mainly taken from the PhD thesis of the second author, 
\cite{Wir04}, achieved under supervision of the first one, and from the
joint preprint \cite{Wir03b}. 

\section{Concepts}
\paragraph{Representation of solutions.}
The Cauchy problem under investigation is invariant under spatial translations. Thus, its 
solution can be represented in terms of Fourier multipliers. If we apply a partial
Fourier transform, we obtain for the solutions to \eqref{eq:CP}
$$ \hat u(t,\xi)=\Phi_1(t,\xi)\hat u_1(\xi)+\Phi_2(t,\xi)\hat u_2(\xi), $$
where the functions $\Phi_i(t,\xi)$ form a fundamental system of solutions to
the ordinary differential equation
\begin{equation}\label{eq:CP_ode}
 \hat u_{tt}+|\xi|^2\hat u+b(t)\hat u_t=0 
\end{equation}
parameterised by the modulus of the frequency variable $\xi$.

\begin{figure}[ht]
\centering
\input{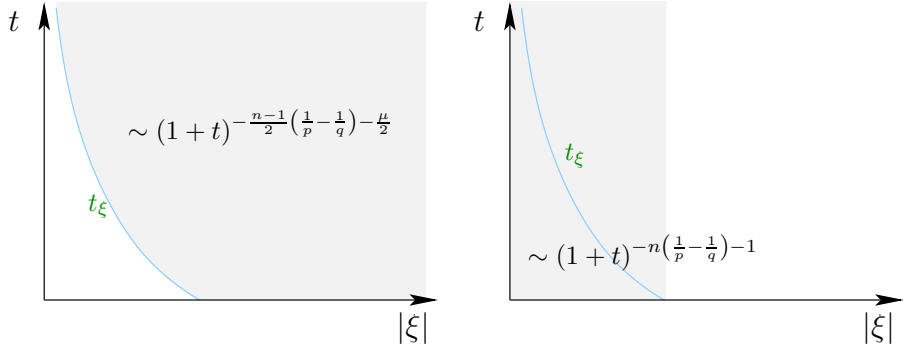}
\caption{Regions of the phase space $\R_+\times\R^n$ determining the $L^p$--$L^q$
decay rate \eqref{eq:Wirth} for small values of $\mu$ (left) and large values of $\mu$ 
(right).}\label{fig:1}  
\end{figure}

\paragraph{Interpretation of known results.}
The main difference between Strichartz' estimate \eqref{eq:Strichartz} and
Matsumura's estimate \eqref{eq:Matsumura} is that the first one is related to properties
of large frequencies and the oscillatory behaviour of the Fourier multipliers, while
the latter one arises from the consideration of small frequencies. 
A similar situation occurs in the estimate \eqref{eq:Wirth}. In \cite{Wir02} the structure of 
the Fourier multipliers $\Phi_i(t,\xi)$ for $b(t)=\frac\mu{1+t}$
is analysed and the resulting decay rates can be related to regions of the
extended phase space $\R_+\times\R^n_\xi$. The change in the $L^p$--$L^q$ decay rate
can be understood as a take-over of the part of the phase space containing the 
``small frequencies''.
In Figure~\ref{fig:1} the different regions are sketched together with
the induced decay rates. 

\paragraph{Effectivity contra non-effectivity.}
For the case of small values of $\mu$, i.e. in the left picture, the dissipation term is
sub-ordinate to the contributions coming from the principal part. We will call a dissipation
term {\em non-effective}, if it does not destroy the basic structure of the Strichartz type 
estimate \eqref{eq:Strichartz} arising from the application of stationary phase method.

In contrast to this, for large values of $\mu$ the influence of the dissipation is much stronger
than the influence of the parameter $|\xi|^2$ in the ordinary differential equation 
\eqref{eq:CP_ode}. In this case there is no need to exploit the oscillatory behaviour of the
Fourier multiplier to get the desired decay rate. If this change of the approach occurs,
we will call the dissipation term {\em effective}.

In order to make the results more precise, we define the {\em energy-operator} 
\begin{equation}\label{eq:EnOp}
 \mathbb E(t) : (\<\D\>u_1,u_2)^T \mapsto (|\D| u(t,\cdot), D_t u(t,\cdot))^T
\end{equation}
and ask for norm estimates of this operator. The operator is normalised in such a way that
it maps $L^2(\R^n,\R^2)\to L^2(\R^n,\R^2)$ for fixed time variable $t$.

\section{Non-effective weak dissipation}
A dissipation term with coefficient function $b=b(t)$ subject to (H1)--(H3) is non-effective if the condition
\begin{itemize}
\item[(NE)] $\displaystyle\limsup_{t\to\infty}b(t)<1$
\end{itemize}
is satisfied. 

\paragraph{The decay estimate.} If we define the auxiliary function 
\begin{equation}\label{eq:lambda_def }
\lambda(t)=\exp\left\{\frac12\int_0^t b(\tau)\d\tau\right\},
\end{equation}
we can describe the asymptotic behaviour of the energy operator $\mathbb E(t)$ and relate it to
the propagator of free waves.

\begin{thm}\label{thm:NE}{\cite[Theorem~3.24]{Wir04}, \cite[Theorem~3.18]{Wir03b}}\\
Assume (H1)--(H3), (NE). Then the energy operator \eqref{eq:EnOp} associated to
the Cauchy problem \eqref{eq:CP} satisfies the norm estimate
\begin{equation}
  ||\mathbb E(t)||_{p,r_p\to q} 
   \lesssim \frac1{\lambda(t)}(1+t)^{-\frac{n-1}2\left(\frac 1p-\frac1q\right)}
\end{equation}
on the conjugate line $pq=p+q$, $q\in[2,\infty]$ and with $r_p>n\left(\frac1p-\frac1q\right)$.
\end{thm}
\begin{proof}[Sketch of the main ideas of the proof] We consider the micro-energy
$U(t,\xi)=(h(t,\xi)\hat u,\D_t\hat u)^T$, where 
$$ h(t,\xi)=\begin{cases} \frac N{1+t},\qquad& (1+t)|\xi|\leq N,\\|\xi|,&(1+t)|\xi|>N,\end{cases}$$
with a suitable constant $N$. This micro-energy satisfies the first order system
$\D_t U=A(t,\xi) U$ and the multiplier $\mathbb E(t,\xi)$ of the energy operator is related to 
the fundamental solution $\mathcal E(t,s,\xi)$
of this system. The main point is to construct a representation of
this fundamental solution, the idea follows \cite{RY00}.
\begin{itemize}
\item We decompose the phase space $\R_+\times\R^n_\xi$ into different zones, a {\em dissipative
zone} containing all $(t,\xi)$ with $|\xi|\leq N b(t)$ and a {\em hyperbolic} one, where
$|\xi|\geq Nb(t)$.
\item In the dissipative zone, we transform this system to an integral equation and prove a uniform
bound for its fundamental solution. This gives under Assumption (NE) the estimate
$$|| \mathcal E(t,s,\xi)||\lesssim \lambda^2(s)/\lambda^2(t). $$
Furthermore, one obtains estimates for derivatives of $\mathcal E(t,s,\xi)$
with respect to $\xi$.
\item In the hyperbolic zone we define symbol classes related to the behaviour of
$b=b(t)$, 
\begin{align*}
  a&(t,\xi)\in S\{m_1,m_2,m_3\}\\
  &\text{iff}\qquad
  \left|\D_t^k\D_\xi^\alpha a(t,\xi)\right|\leq C_{k,\alpha} |\xi|^{m_1-|\alpha|} b(t)^{m_2}
  \left(\frac1{1+t}\right)^{m_3+k}
\end{align*}
for all $k\in\N$, $\alpha\in\N^n$ and all $(t,\xi)$ inside the hyperbolic zone.
Thus $b(t)\in S\{0,1,0\}$ by (H3). The symbol classes satisfy natural rules of symbolic
calculus.
\item
We apply several steps of diagonalization to this system with respect to the symbol hierarchy
$S\{-k,1,k\}\hookrightarrow S\{0,1,0\}$, $k\geq0$. 
This 
gives for each $k\in\N$ and a suitable choice of the zone constant $N$ an invertible matrix
$N_k(t,\xi)\in S\{0,0,0\}$, such that 
$$\mathcal E_k(t,s,\xi)=N_k^{-1}(t,\xi)\mathcal E(t,s,\xi)N_k(s,\xi) $$
satisfies inside the hyperbolic zone 
$$ \D_t\mathcal E_k-\big(\mathcal D(\xi)+F(t,\xi)+R_k(t,\xi)\big)\mathcal E_k=0 $$
with $\mathcal D(\xi)=\mathrm{diag}\,(-|\xi|,|\xi|)$, $F(t,\xi)$ diagonal
with $F(t,\xi)-\frac i2 b(t)I\in S\{-1,1,1\}$ and $R_k(t,\xi)\in S\{-k,1,k\}$.
\item
Now, we can use the Peano-Baker formula to represent the fundamental solution of the 
transformed system and also to estimate a finite number of derivatives with respect to the 
frequency variable $\xi$. 
\end{itemize}
We can use the obtained representation of $\mathcal E(t,s,\xi)$ and, therefore,
that of $\mathbb E(t,\xi)$ to deduce the desired decay estimate. For this, we use
the stationary phase method in combination with Marcinkiewicz multiplier theorem.
\end{proof}

\paragraph{Sharpness.}
The representation of the multiplier $\mathbb E(t,\xi)$ allows us to obtain even more. If we denote
by $\mathbb E_0(t)$ the unitary propagator of free waves in the energy space, i.e.
\begin{equation}
 \mathbb E_0(t) : (|\D|\tilde u(0,\cdot),\D_t \tilde u(0,\cdot))^T\mapsto
   (|\D| \tilde u(t,\cdot), D_t \tilde u(t,\cdot))^T
\end{equation}
for a solution $\tilde u=\tilde u(t,x)$ to the free wave equation $u_{tt}-\Delta u=0$,
we obtain an asymptotic equivalence of $\lambda(t)\mathbb E(t)$ and $\mathbb E_0(t)$ in the 
following sense.

\begin{thm}{\cite[Theorem 3.26]{Wir04}}\\
Assume (H1)--(H3), (NE). Then the limit 
$$ W_+ = \slim_{t\to\infty} \lambda(t)(\mathbb E_0(t))^{-1}\mathbb E(t) $$
exists as strong limit in $L^2\to L^2$ and defines a bounded and injective translation invariant
operator with dense range.
\end{thm}

This result is closely related to the construction of the M\o{}ller wave operator in scattering 
problems, see e.g. \cite{LP73} or \cite{Mel95} and the discussion in Example~\ref{expl3.1}. The operator
$W_+$ associates to Cauchy data $(\<\D\>u_1,u_2)^T$ of the damped wave equation \eqref{eq:CP} data
$(|\D|\tilde u_1,\tilde u_2)^T$ to the free problem, such that the modified solution $\lambda(t)u(t,x)$
of the damped problem  and the free solution $\tilde u(t,x)$ are asymptotically equivalent.

Thus this theorem may be used to obtain the sharpness of the above given energy estimate. Especially it
provides us with a lower bound for the $L^2$--$L^2$ decay rate. 

\begin{cor} Assume (H1)--(H3), (NE). Then it holds
$$ ||\mathbb E(t)||_{2\to2}\sim \frac1{\lambda(t)}. $$
\end{cor}

We conclude this section with several examples to underline the previous statements. 
\begin{expl}\label{expl3.1}
If we assume that $b(t)\in L^1(\R_+)$, we have $\lambda(t)\sim 1$ and the estimates simplify
to the known Strichartz' decay estimates for free waves. We obtain a scattering result in 
the energy space, which is the counterpart to the result of Mochizuki, \cite{Moc76}, \cite{MN96},
and relates the non-decay to zero of the energy to the conservation of energy for free waves.
\end{expl}
\begin{expl}\label{expl:ItLogNE}
If we set 
$$ b(t)=\frac\mu{(1+t)\ln(e+t)\cdots\ln^{[m]}(e^{[m]}+t)} $$
with iterated logarithms $\ln^{[k+1]}(t)=\ln(\ln^{[k]}(t))$ and
exponentials \mbox{$e^{[k+1]}=e^{(e^{[k]})}$}, we can obtain arbitrarily small decay rates for the 
energy. We have
$$ \lambda(t)\sim \big(\ln^{[m]}(e^{[m]}+t)\big)^{\frac\mu2} $$
and therefore,
$$ ||\mathbb E(t)||_{p,r_p\to q}\lesssim \big(\ln^{[m]}(e^{[m]}+t)\big)^{-\frac\mu2}(1+t)^{-\frac{n-1}2\left(\frac1p-\frac1q\right)} $$
and 
$$ ||\mathbb E(t)||_{2\to2}\sim\big(\ln^{[m]}(e^{[m]}+t)\big)^{-\frac\mu2}.$$
\end{expl}
\begin{expl}
If we set $b(t)=\frac\mu{1+t}$ with $\mu\in(0,1)$, we obtain the estimate \eqref{eq:Wirth} for 
this case, i.e.
$$ ||\mathbb E(t)||_{p,r_p\to q}
\lesssim (1+t)^{-\frac{n-1}2\left(\frac1p-\frac1q\right)-\frac\mu2}, $$
together with the description of the energy decay
$$ ||\mathbb E(t)||_{2\to2}\sim (1+t)^{-\frac\mu2}. $$
The value $\mu=1$ is exceptional for estimates of the solution itself, see \cite{Wir02}. 
Furthermore, following \cite[Example 3.4]{Wir04}, we obtain the above given estimates also for 
$\mu\in(1,2]$.  
\end{expl}

\section{Effective dissipation}
Now we devote our study to dissipation terms which lie above the intermediate case $\mu/(1+t)$. 
We assume, that the condition
\begin{itemize}
\item[(E)] $tb(t)\to\infty$ as $t\to\infty$
\end{itemize}
is satisfied.

\paragraph{Transformation of the problem.} 
In this case we have to look more carefully to the behaviour of small frequencies. Following the
treatment in \cite{Rei00}, \cite{Mat76}, the basic idea is a transformation of the dissipative problem
to a Klein-Gordon type equation. This can be achieved by the consideration of the new function
\begin{equation}\label{eq:v_def}
 v(t,x)=\lambda(t)u(t,x),
\end{equation}
such that
\begin{equation}\textstyle
 v_{tt}-\Delta v=\big( \frac14 b^2(t)+\frac12 b'(t)\big)v.
\end{equation}
Under the assumptions of non-effective weak dissipation, (H1)--(H3), (NE), 
$-b'(t)\lesssim b(t)/(1+t)$ dominates $b^2(t)$ and the potential term is positive for large time 
$t$. In our case $b^2(t)$ dominates $-b'(t)$ and at least for small frequencies we feel a
negative potential. 

\paragraph{The decay estimate.} For effective dissipation terms the structure of the estimate
changes completely. It holds:

\begin{thm}\label{thm:Eff}{\cite[Theorem 4.25]{Wir04},\cite[Theorem~4.8]{Wir03b}}\\
Assume (H1)--(H3), (E). Then the energy operator \eqref{eq:EnOp} associated to
the Cauchy problem \eqref{eq:CP} satisfies the norm estimate
\begin{equation}
  ||\mathbb E(t)||_{p,r_p\to q} 
   \lesssim \left(1+\int_0^t\frac{\d\tau}{b(\tau)}\right)^{-\frac{n}2\left(\frac 1p-\frac1q\right)-\frac12}
\end{equation}
on the conjugate line $pq=p+q$, $q\in[2,\infty]$ and with $r_p>n\left(\frac1p-\frac1q\right)$.
\end{thm}
\begin{figure}
\centering\input{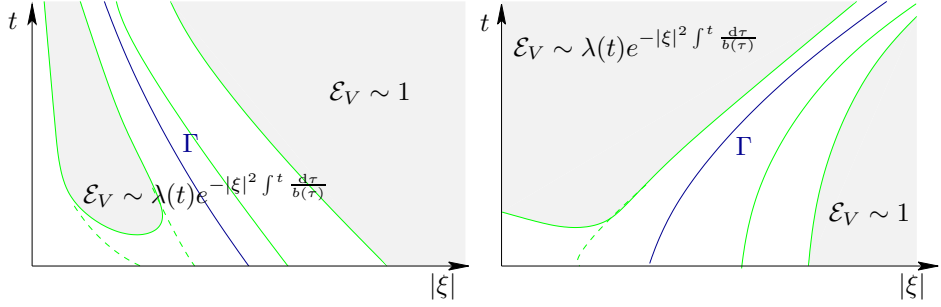}
\caption{Parts and zones used in the case of effective dissipation, left
for decreasing $b=b(t)$ and on the right for increasing dissipation.}\label{fig:Eff}
\end{figure}
\begin{proof}[Sketch of the main ideas of the proof]
Like in the proof of Theorem~\ref{thm:NE} we construct a representation of the fundamental solution
$\mathcal E_V(t,s,\xi)$ after applying transformation \eqref{eq:v_def}. The approach
follows basically \cite{Rei00}, we summarise the main steps.
\begin{itemize}
\item
 We decompose the extended phase space $\R_+\times\R^n$ into the {\em elliptic part}, where
$m(t,\xi)=|\xi|^2-\frac14b^2(t)$ is positive, and the {\em hyperbolic part}, where it is negative. 
In both parts we transform the problem into system form using
$V(t,\xi)=\big(\sqrt{|m(t,\xi)|}\hat v(t,\xi),\D_t\hat v(t,\xi)\big)^T$.
The curve $\Gamma=\{2|\xi|=b(t)\}$ we denote as {\em separating curve}, $\sqrt{|m(t,\xi)|}$
is a measure of the distance of $(t,\xi)$ to the separating curve $\Gamma$.
\item
In both parts we introduce zones, a {\em hyperbolic} and an {\em elliptic zone}, to stay away 
from the separating curve. Furthermore, we introduce symbol classes. They are defined as
\begin{align*}
  a&(t,\xi)\in S_{ell}\{m_1,m_2,m_3\}, \qquad \big(\in S_{hyp}\{m_1,m_2,m_3\}\big), \\
  &\text{iff}\qquad
  \left|\D_t^k\D_\xi^\alpha a(t,\xi)\right|\leq C_{k,\alpha} 
  \left(\sqrt{|m(t,\xi)|}\right)^{m_1-|\alpha|} (b(t))^{m_2}
  \left(\frac1{1+t}\right)^{m_3+k}
\end{align*}
for all $k\in\N$, $\alpha\in\N^n$ and all $(t,\xi)$ inside the elliptic (hyperbolic) zone.
Again these symbol classes satisfy natural rules of symbolic calculus. Besides these two zones
we introduce a {\em reduced zone} in the neighbourhood of the separating curve and 
near the $t$-axis
we can apply ideas as in the {\em dissipative zone} of the proof of 
Theorem~\ref{thm:NE}.
\item Inside the hyperbolic zone we show that $||\mathcal E_V(t,s,\xi)||\sim 1$. More structural
properties of $\mathcal E_V(t,s,\xi)$ may be obtained by a diagonalization procedure like it is
used in the proof of Theorem~\ref{thm:NE}.
\item Inside the elliptic zone we apply a similar diagonalization procedure to decouple the system
modulo $S_{ell}\{-1,0,2\}$. We can not use Peano-Baker formula to estimate the fundamental solution
$\mathcal E_{V,1}(t,s,\xi)$
of this transformed system. The basic idea is to rewrite the system as an integral equation for 
$\mathcal Q_{ell,1}(t,s,\xi)$ from the approach
$$ \mathcal E_{V,1}(t,s,\xi)=\frac{\sqrt{|m(t,\xi)|}}{\sqrt{|m(s,\xi)|}} 
\exp\left\{\int_s^t\sqrt{|m(\tau,\xi)|}\d\tau\right\} \mathcal Q_{ell,1}(t,s,\xi) $$
and to prove uniform bounds for $\mathcal Q_{ell,1}(t,s,\xi)$, $t\geq s$, in the elliptic zone.
\item
The reduced zone is of minor influence, so a rough estimate is sufficient.  Transformation
back to $\hat u(t,\xi)$ yields in the hyperbolic zone a decay like $\lambda^{-1}(t)$ and
cancels the increasing behaviour of the exponential from the elliptic zone partly. It holds
$$
\exp\left\{\int_s^t\sqrt{|m(\tau,\xi)|}-\frac12 b(\tau)\d\tau\right\}
\leq \exp\left\{-|\xi|^2\int_s^t\frac{\d\tau}{b(\tau)}\right\}.
$$
\item
The decay rate is determined from the elliptic part. There is no need to apply the 
stationary phase method, the $L^1$--$L^\infty$ decay rate can be estimated by the $L^1$-norm
of the multiplier $\mathbb E(t,\xi)$ over the elliptic part.
\end{itemize}
\end{proof}

We proceed with two examples to give an impression of the obtained decay rates. The $L^2$--$L^2$
rate of the first one is known from \cite{Ues79}, the $L^p$--$L^q$ rates extent
the well-known Matsumura estimate.
\begin{expl}
If we consider $b(t)=(1+t)^\kappa$ with $\kappa\in(-1,1)$, we obtain the estimate
$$
 ||\mathbb E(t)||_{p,r_p\to q}\lesssim (1+t)^{(\kappa-1)\left(\frac n2\left(\frac1p-\frac1q\right)+\frac12\right)}.$$
This estimate fits to \eqref{eq:Wirth} for $\kappa\to-1$ and contains also the estimate
of Matsumura, \eqref{eq:Matsumura}, as a special case.
\end{expl}
\begin{expl}
If we set $b(t)=1+t$, we obtain the logarithmic estimate
$$  ||\mathbb E(t)||_{p,r_p\to q}\lesssim \big(\log(e+t)\big)^{-\frac n2\left(\frac1p-\frac1q\right)-\frac12}.$$
A construction related to Example~\ref{expl:ItLogNE} leads to decay rates of arbitrarily small
logarithmic order.
\end{expl}

\paragraph{Sharpness.}
The basic idea of the proof of Theorem~\ref{thm:Eff} to construct the leading terms of the
representation of solutions hints that the achieved estimates are indeed sharp. The questions 
related to this sharpness will conclude this section on effective dissipation terms. We start 
with the following special case, called the case of {\em over-damping}.

\begin{expl}
 If we assume $1/b(t)\in L^1(\R_+)$  Theorem~\ref{thm:Eff} trivialises to the 
 (in view of condition (H1) obvious) estimate $ ||\mathbb E(t)||_{p,r_p\to q} \lesssim 1$.
\end{expl}

From the representation of the operator $\mathbb E(t)$ as Fourier multiplier we can 
conclude that this estimate is indeed sharp. 

\begin{thm}{\cite[Theorem 4.27 and 4.31]{Wir04}}\\
Assume (H1)--(H3), (E) and  $1/b(t)\in L^1(\R_+)$. Then for $u_1\in H^1(\R^n)$ and
$u_2\in L^2(\R^n)$ the solution $u(t,x)$ to \eqref{eq:CP} converges in $H^1(\R^n)$
to the asymptotic state
$$ u(\infty,x)=\lim_{t\to\infty} u(t,x), $$
which is a real-analytic function in $x$ and non-zero for non-zero initial data.
\end{thm}

As a consequence we see that in the case of over-damping the solution and its spatial
derivatives can not decay to zero in $L^2(\R^n)$ and, therefore, also not (locally)
in $L^\infty(\R^n)$. 

For the remaining effective dissipation terms the question of sharpness is more involved.
Following \cite{HN03} for the case $b(t)=\frac\mu{1+t}$, $\mu>2$ 
or  $b(t)=1$, we see that the hyperbolic energy
\begin{equation}
 E(u;t)=\frac12\int_{\R^n} \big(|\nabla u|^2+|u_t|^2\big)\d x 
\end{equation}
decays slightly faster than the estimate of Theorem~\ref{thm:Eff} predicts. It holds
\begin{align*}
 &E(u;t)=o(t^{-2}),\quad t\to\infty, \qquad\qquad && b(t)=\frac\mu{1+t},\quad\mu>2,\\
 &E(u;t)=o(t^{-1}),\quad t\to\infty, \qquad\qquad && b(t)=1.
\end{align*}
The $L^2$--$L^2$ norm-estimate of Theorem~\ref{thm:Eff} and these estimates 
in the strong topology, i.e. in dependence of particular data, are both sharp. 
This is a consequence of the representations of the 
multiplier inside the elliptic part.

\begin{thm}\label{thm:6}{\cite[Theorem 5.1 and Corollary 4.12]{Wir04}}\\
 Assume (H1)--(H3), (E) and $1/b(t)\not\in L^1(\R_+)$. Then
$$
  \slim_{t\to\infty} \sqrt{1+\int_0^t\frac{\d\tau}{b(\tau)}}\;\mathbb E(t)=0
$$
in $L^2(\R^n)$ and there exists no function $\omega(t)$ with $\omega(t)\to0$
as $t\to\infty$, such that
$$ \sqrt{1+\int_0^t\frac{\d\tau}{b(\tau)}}\; ||\mathbb E(t)||_{2\to2} \leq \omega(t). $$
\end{thm}
\begin{proof}[Sketch of the main ideas of the proof]
The proof is based on two facts.
\begin{itemize}
\item The representation of $\mathcal Q_{ell,1}(t,s,\xi)$ from the proof of
Theorem~\ref{thm:Eff} implies that this matrix tends to a nonzero limit as $t\to\infty$
locally uniform in $\xi$. This can be used to deduce
$$ ||\mathbb E(t)||_{2\to2}
=||\mathbb E(t,\cdot)||_\infty\sim \left(1+\int_0^t\frac{\d\tau}{b(\tau)}\right)^{-\frac12} $$
and therefore the norm-estimate is sharp.
\item On the other hand, if $1/b(t)\not\in L^1(\R_+)$ we obtain for data from the
subspace $V_c=\{u\in L^2\;|\; \mathrm{dist}\,(0,\mathrm{supp}\, u)\geq c \}$ a stronger decay rate.
Because the union $M=\bigcup V_c$ is dense in $L^2$, we can apply the theorem of Banach-Steinhaus
to obtain the strong convergence.
\end{itemize}
\end{proof}

\section{Concluding remarks}
As the basic idea to obtain the collected $L^p$--$L^q$ decay estimates for the energy of the solution to
the damped problem \eqref{eq:CP} we used a precise construction of the main terms of the representation
of solutions. These constructions may also be used to handle several related problems. To conclude this expository
article we give some remarks on such applications; the depiction cannot be regarded as a complete one.

\paragraph{Estimates for the solution itself.} The definition of the micro-energies $U(t,\xi)$ in the proof of Theorem~\ref{thm:NE}
and $V(t,\xi)$ in the one of Theorem~\ref{thm:Eff} can be used to extract the Fourier transform of the solution 
$\hat u(t,\xi)$ and, thus, to deduce also estimates for it. Such estimates are given in \cite{Wir02} for the case
of $b(t)=\mu/(1+t)$ and in \cite[Chapter 5.2]{Wir04} in generality.

If we define the {\em solution operator}
\begin{equation}\label{eq:SolOp}
 \mathbb S(t) : (u_1,\<\D\>^{-1}u_2)^T \mapsto u(t,\cdot),
\end{equation}
we can formulate the following two theorems. In the case of non-effective dissipation
the estimates are closely related to the corresponding estimates for free waves.

\begin{thm}{\cite[Theorem~5.9]{Wir04}}\\
Assume (H1)--(H3) together with (NE). Then the $L^p$--$L^q$ estimate
$$ ||\mathbb S(t)||_{p,r_p\to q} \lesssim
\begin{cases} 
\frac1{\lambda(t)}(1+t)^{-\frac{n-1}2\left(\frac1p-\frac1q\right)},\qquad&p<p^*\\
\frac1{\lambda^2(t)}(1+t)^{1-n\left(\frac1p-\frac1q\right)},&p\geq p^* \end{cases}$$
holds for dual indices $q\in[2,\infty]$, $pq=p+q$ and with $r_p>n\left(1/p-1/q\right)$. The critical value $p^*$
is chosen from $(n+1)(1/p^*-1/2)= \liminf_{t\to\infty}(1-\log_t\lambda(t))$.
\end{thm}

For $p$ and $q$ near $2$ the estimate is determined from the dissipative zone in opposite to the energy estimate
in this case. As consequence we obtain similar to Theorem~\ref{thm:6}
that $||u(t,\cdot)||_{2}=o(t/\lambda^2(t))$ for fixed initial data. 

For effective dissipation the structure of the estimate is related to the one of Theorem~\ref{thm:Eff}.

\begin{thm}{\cite[Theorem~5.11]{Wir04}}\\
Assume (H1)--(H3) together with (E). Then  the $L^p$--$L^q$ estimate
$$ ||\mathbb S(t)||_{p,r_p\to q} \lesssim
 \left(1+\int_0^t\frac{\d\tau}{b(\tau)}\right)^{-\frac n2\left(\frac1p-\frac1q\right)} $$
holds for dual indices $q\in[2,\infty]$, $pq=p+q$ and with $r_p>n\left(1/p-1/q\right)$.
\end{thm}

\paragraph{Diffusive structure.} In the case of effective dissipation we obtained in the elliptic part as main 
term in the representation of solutions the expression
$$ \exp\left\{-|\xi|^2\int_0^t\frac{\d\tau}{b(\tau)}\right\},$$
which turns out to be the Fourier multiplier representing the solution of the associated {\em parabolic problem}
\begin{equation}\label{eq:parCP}
   b(t)w_t=\Delta w,\qquad w(0,\cdot)=w_0.
\end{equation}
It is a natural question under which assumptions on the coefficient function the solutions to the hyperbolic 
problem \eqref{eq:CP} and to this parabolic surrogate \eqref{eq:parCP} are asymptotically equivalent. For $b(t)=1$
this relation was treated in \cite{Nis97}, \cite{HM00}, \cite{Nar04} and several other papers and is referred to as 
the {\em diffusion phenomenon}. Using the constructed representation of solutions the diffusion phenomenon can be 
extended to a neighbourhood of the case $b(t)=1$, \cite[Chapter 5.4]{Wir04}.

Related to this diffusive structure is the content of Theorem~\ref{thm:6}. Using further assumptions on the data
which are effective near the exceptional frequency $\xi=0$ we can obtain improved decay rates of the energy under
the condition $1/b(t)\not\in L^1(\R_+)$, \cite[Chapter 5.1]{Wir04}. Examples for such assumptions are
\begin{itemize}
\item data from $H^s\cap L^p$ with $p\in[1,2)$ like in \cite{Mat76},
\item data satisfying weight conditions like in \cite{01931006}. 
\end{itemize}

\paragraph{Estimates for higher order energies.}
In the case of free waves we can differentiate the equation with respect to all variables;
so energies of higher order are preserved like the usual first order energy. On the other hand, 
for the damped wave equation the results of \cite{Mat76} give stronger decay rates for higher order derivatives. 

Our representation of solutions can also be used to deduce also such estimates. The situation occurs as follows,
\cite[Chapter 5.3]{Wir04}:
\begin{itemize}
\item in the case of non-effective dissipation, higher order derivatives behave like first order derivatives
and the same estimates are valid, i.e. we have under the assumptions (H1)--(H3), (NE) the norm estimate
$$ ||\D_t^k\D_x^\alpha u(t,\cdot)||_q \lesssim \frac1{\lambda(t)} (1+t)^{-\frac{n-1}2\left(\frac1p-\frac1q\right)}
 ||(u_0,\<\D\>^{-1} u_1)||_{L^{p,r_p+k+|\alpha|}} $$
for $k+|\alpha|\geq1$ and with $q\in[2,\infty]$, $pq=p+q$ and $r_p>n(1/p-1/q)$,
\item in the case of effective dissipation the decay order depends on the number of spatial and time derivatives,
where time-derivatives bring more improvement than spatial ones. As special case we refer to the estimates of
Matsumura, \cite{Mat76}, which imply under the same assumptions as above
$$ ||\D_t^k\D_x^\alpha u(t,\cdot)||_q\lesssim
   (1+t)^{-\frac n2\left(\frac1p-\frac1q\right)-k-\frac{|\alpha|}2} 
   ||(u_0,\<\D\>^{-1} u_1)||_{L^{p,r_p+k+|\alpha|}}. $$
\end{itemize}

\paragraph{Applications to nonlinear problems.}
$L^p$--$L^q$ decay estimates are a classical tool to treat nonlinear problems. For the special case
\begin{equation}
  u_{tt}-\Delta u+\frac{\mu}{1+t}u_t=f(u'),\qquad u'=(u_t,\nabla u)^T
\end{equation}
with a nonlinearity $f\in C^\infty(\R^{n+1})$, $f(0)=0$, $\D^\alpha f(0)=0$ for $|\alpha|=1$ and 
for effective dissipation, $\mu>2$, this can be done using standard arguments. In this case it is possible to show that
small data solutions exist globally without further assumptions on the nonlinearity. 

Thus we can treat arbitrary quadratic nonlinearities if we allow this kind of dissipation terms. If we consider
nonlinear perturbations of the free wave equation the situation is quite different. Then for space dimensions
$n=2$ or $n=3$ we need further conditions on the structure of the nonlinearity like Klainerman's famous 
Null condition, \cite{Kla83}. The John example $u_{tt}-\Delta u=f(u')=(u_t)^2$ leads to a blow-up
for arbitrary non-zero small data solutions, \cite{Joh81}. 

\paragraph{Monotonicity of the coefficient function} is a technical assumption in order to include
it in the estimates for the symbolic calculus. In the case of non-effective dissipation we can replace (H2) and (H3) by (H3') $$|b^{(k)}(t)| \leq C_k(1+t)^{-k-1}.$$ 

In the case of effective dissipation the function $b(t)$ is used in the definition of the separating
curve between the hyperbolic and the elliptic part and, therefore, the monotonicity is used in an essential
way. To treat also non-monotonous coefficients we need a 
monotonous comparison function $\gamma(t)$ subject to (H1)--(H3), (E) and assume for the coefficient function
$b(t)$ the conditions (H1), (H3) (with $b(t)$ replaced by $\gamma(t)$) and (H$\gamma$) 
$$ |b(t)-\gamma(t)|\leq (1+t)^{-1}.$$

These more general assumptions on the coefficient $b=b(t)$ in combination with the corresponding symbol classes 
allow us to apply the sketched diagonalization procedure and to derive corresponding expressions of the main terms
of the representation of solutions. In this overview article we preferred to use the simpler assumptions (H1)--(H3) 
instead in order to emphasise the general philosophy behind the results. For the general treatment we refer to
the PhD thesis \cite{Wir04} and a planned series of forthcoming papers on the subject.
\bibliographystyle{abbrv}
\bibliography{../database.bib}
\end{document}